\documentclass{article}
\usepackage[utf8]{inputenc}
\usepackage{amsmath}
\usepackage{amssymb}
\usepackage{amsthm}
\usepackage{mathtools}
\usepackage{tikz-cd}
\usepackage[hyphens]{url}
\usepackage{hyperref}
\usepackage[hyphenbreaks]{breakurl}
\usepackage[style=numeric,sorting=nty]{biblatex}
\newtheorem{theorem}{Theorem}
\newtheorem{lemma}{Lemma}
\newtheorem{proposition}{Proposition}
\newtheorem{definition}{Definition}

\title{Computing Tangent Spaces to Eigenvarieties}
\author{James Rawson\thanks{The author is supported by the Warwick Mathematics Institute Centre for Doctoral Training, and gratefully acknowledges funding from the UK Engineering and Physical Sciences Research Council (Grant number: EP/W523793/1)}}

\begin{document}

\maketitle

\begin{abstract}
We develop an effective algorithm to compute the derivative of a Bianchi modular form with respect to weight space as it varies in a $p$-adic family. This method is entirely local at the modular form, and does not compute the family anywhere outside an infinitesimal neighbourhood. We numerically verify some conjectures surrouding smoothness of the eigenvariety (equivalently, uniqueness of families) and the ``direction over weight space'' of the family. The methods are also applied to study elliptic modular forms and their $\mathcal{L}$-invariants.
\end{abstract}

\section{Introduction}
By work of Hida and Coleman, it is known that modular forms vary in (one-dimensional) $p$-adic families, that is, there exist $p$-adic analytic functions $a_i(k)$ such that $\sum_{n \geq 0} a_n(k) q^n$ is a modular form of weight $k$ for ``nice enough'' $k$. These can be given a geometric interpretation via the eigenvariety construction. Many properties of these $p$-adic families are well-understood, for example, a family passing through a given eigenform (satisfying some mild conditions) is unique (equivalently, the eigenvariety is smooth at the point corresponding to the eigenform). Similarly, Bianchi cusp forms vary in one-dimensional $p$-adic families, but a lot less is known about their behaviour. As Bianchi modular forms are described by two weights, one question is how do these one-dimensional families vary over the weight space? Calegari and Mazur conjecture that classical families (families containing infinitely many classical forms) should arise only for forms coming from base change or from complex multiplication \cite{calmaz}. As cusp forms in the Bianchi setting are constrained to parallel weight, this is equivalent to families through forms that are not CM or base change are not contained in the parallel weight direction. This would certainly be implied by the derivatives of the weights of the $p$-adic family not being equal. A second question is whether Bianchi modular forms belong to unique families. As well as having inherent interest for the theory of Bianchi modular forms, this is hypothesis is an assumption in a resolution of the analytic rank 0 case of the Birch and Swinnerton-Dyer conjecture for modular elliptic curves over quadratic imaginary fields \cite{loezer}.

Throughout this paper, $K$ is a quadratic imaginary number field, $p$ is a rational prime that splits in $K$ and $\Omega$ is an adelic congruence subgroup contained in $\Gamma_0(p)$.
\subsection{Results}
In this paper, we study tangent spaces to eigenvarieties in terms of cohomology groups with coefficients given by modules over the dual numbers. From this approach, we were able to recover a variant of a result of Hansen and Newton on the dimensions of irreducible components of the Bianchi eigenvariety \cite{hannew} (Proposition B.1). Namely, we prove the following.
\begin{theorem}
Let $\Phi$ be a Bianchi cusp form with sub-critical slope, then the dimension of the tangent space to the eigenvariety at $\Phi$ is at least 1 dimensional.
\end{theorem}
The methods used to prove this can be made completely explicit, which gives an algorithm to compute the dimension of the tangent space to the eigenvariety at a given point. When the computed bound is 1, it is also able to compute its image in the tangent space to the weight space (the ``deformation direction'' of the form). This can all be done in polynomial time We implemented this algorithm, and ran it on a few examples. In each of these cases, the output of the program was that the tangent space is 1-dimensional, and the $p$-adic family through the form does not contain the diagonal, verifying the conjectures discussed before.

The methods are not restricted to the Bianchi eigenvariety, as we also apply them to the Coleman-Mazur eigencurve. In this setting, we can compute more detailed information about the tangent spaces. In particular, we compute the derivatives of the Hecke eigenvalues in the $p$-adic family, and consequently, compute $\mathcal{L}$-invariants. We illustrate this with a few examples, including a form of weight 4.

\subsection{Structure of the Paper}
Section 1 is a review of standard definitions and conventions surrounding overconvergent Bianchi modular forms and Bianchi modular symbols. Section 2 provides the theoretical background for the algorithm and explores certain elementary properties of tangent spaces to the Bianchi eigenvariety. Section 3 outlines the algorithm and its complexity along with its output in a few cases which verifies both the conjecture around smoothness and the conjecture of Calegari and Mazur. In Section 4, we apply the methods of Section 2 and 3 to elliptic modular forms, these are used to compute $\mathcal{L}$-invariants of modular forms as well as derivatives of Hecke eigenvalues with respect to the weight of the form.

\section{Overconvergent Bianchi Modular Symbols}
We start by defining a pair of $\mathrm{SL}_2(K)$ modules, following the coventions used in \cite{williams}.

\begin{definition}
Let $\mathcal{P}_{k, k}(R)$ denote polynomials in $x$ and $y$ of degree at most $k$ in each variable with coefficients in a $K$-algebra, $R$. This comes with the following left action.
$$\begin{pmatrix}a & b \\ c & d\end{pmatrix}\cdot f(x, y) = (a + cx)^k (\bar{a} + \bar{c}y)^k f\left(\frac{dx + b}{cx + a}, \frac{\bar{d}y + \bar{b}}{\bar{c}x + \bar{a}}\right)$$
The module $V_{k, k}(R)$ is the dual of $\mathcal{P}_{k, k}(R)$, with a right action given in terms of this left action.
$$\left(\mu \cdot \gamma\right)(f) = \mu(\gamma \cdot f)$$
\end{definition}

With these definitions, we can compare Bianchi modular forms to cohomology.
\begin{theorem}[Eichler-Shimura-Harder Isomorphism]
There is an isomorphism between the space of Bianchi cusp forms of weight $(k + 2, k + 2)$, $S_{k, k}(\Omega)$, and $H^1_{c, \mathrm{cusp}}(Y_{\Omega}, V_{k, k})$, where $Y_{\Omega}$ is the following quotient.
$$\mathrm{GL}_2(K) \backslash \mathrm{GL}_2(\mathbb{A}_K) / (\Omega\; \mathrm{SU}_2(\mathbb{C})\; \mathbb{C}^{\times})$$
Moreover, this isomorphism respects the action of the Hecke operators.
\end{theorem}

For brevity, we will suppress $Y_{\Omega}$ in cohomology groups. For use in later explicit calculations, we make the following definition.
\begin{definition}
For any discrete subgroup $\Gamma \leq \mathrm{SL}_2(K)$, modular symbols of level $\Gamma$ taking values in a right $\Gamma$-module, $V$, are elements of $\mathrm{Symb}_{\Gamma}(V) := \mathrm{Hom}_{\Gamma}(\mathrm{Div}^0(\mathbb{P}^1(K)), V))$, with $\Gamma$-equivariance given by $\phi(\gamma D) = \phi(D) \cdot \gamma^{-1}$.
\end{definition}
These are isomorphic to the cohomology groups $H^1_{c}(\Gamma \backslash \mathbb{H}_3, V)$, where $\mathbb{H}_3$ is hyperbolic 3-space. The group $H^1_c(Y_{\Omega}, V)$ is the direct sum $\oplus_{i = 1}^{h_K} H^1_c(\Gamma_i \backslash \mathbb{H}_3, V)$, where $\Gamma_i = \mathrm{SL}_2(K) g_i \Omega g_i^{-1}$, and the $g_i$ come from representatives for the class group. The sum of the modular symbols spaces will be denoted $\mathrm{Symb}_{\Omega}(V)$.

For the purposes of $p$-adic variation, we introduce two additional modules.
\begin{definition}
Let $\Sigma_0(p) = \{ \begin{psmallmatrix}a & b \\ c & d \end{psmallmatrix} \in M_2(\mathbb{Z}_p) \mid a \in \mathbb{Z}_p^{\times}, c \in p\mathbb{Z}_p, ad - bc \neq 0\}$.

The space of $R$-valued analytic functions (for $R$ a Banach $\mathbb{Q}_p$-algebra), $\mathcal{A}(R)$, on $\mathbb{Z}_p^2$ is a module for $\Sigma_0(p)^2$. For any pair of characters $\chi_1, \chi_2 : \mathbb{Z}_p^{\times} \mapsto R^{\times}$, the action of $\Sigma_0(p)^2$ is given by $((\sigma_1, \sigma_2) \cdot f)(x, y) = \chi_1(a_1 + c_1x) \chi_2(a_2 + c_2 y) f\left(\frac{d_1 x + b_1}{a_1 + c_1 x}, \frac{d_2 y + b_2}{a_2 + c_2 y}\right)$. 

Similarly to before, there is also the dual space, $\mathbb{D}(R) = \mathrm{Hom}(\mathcal{A}(\mathbb{Q}_p), \mathbb{Q}_p) \hat{\otimes} R$. This is equipped with the natural action. When the weight needs to be made explicit, this will be denoted by $\mathbb{D}_{\chi_1, \chi_2}$.

The group $\Gamma_0(p)$ includes into $\Sigma_0(p)^2$ via $\sigma \mapsto (\sigma, \bar{\sigma})$.
\end{definition}
Modular symbols valued in such modules are called overconvergent modular symbols. The Hecke algebra, $\mathbb{T}$, acting on these overconvergent modular symbols will be generated by $T_{\mathfrak{q}}$ (for $\mathfrak{q}$ not dividing the level), $U_{\mathfrak{q}}$ ($\mathfrak{q}$ dividing the level), and the symbols $\langle a \rangle = \begin{psmallmatrix}a & 0 \\ 0 & a \end{psmallmatrix}$, for all $a \in K^{\times}$. The last family of operators function as proxies for the weight, since these act on distributions via $\mu \mapsto \chi_1(a) \chi_2(\bar a) \mu$.

The final definition is that of the eigenvariety, following Bella\"iche's Eigenbook \cite{eigen}.
\begin{definition}
Weight space, $\mathcal{W}$, is the rigid analytic space given by the following functor on Banach algebras over $\mathbb{Q}_p$: $\mathcal{W}(R) = \mathrm{Hom}(\mathbb{Z}_p^{\times}, R^{\times})$. The Bianchi eigenvariety is given as a covering of $\mathcal{W}^2$, where above $(\chi_1, \chi_2) \in \mathcal{W}(R)^2$, the eigenvariety is given by the finite slope part of the spectrum of $\mathbb{T}$ acting on $H^1_c(\mathbb{D}_{\chi_1, \chi_2}(R))$. 
\end{definition}

\section{Tangent Spaces to the Eigenvariety}
To compute the tangent space to a rigid space, $\mathcal{X}$, at a point defined over $L$, we study maps $\mathrm{Spm}(L[\varepsilon]/(\varepsilon^2)) \to \mathcal{X}$ such that sending $\varepsilon$ to 0 gives back the inclusion of $\mathrm{Spm}(L)$. More precisely, the dimension of the tangent space at a point is the dimension of this space of maps. The ring $L[\varepsilon]/(\varepsilon^2)$ will be abbreviated to just $L[\varepsilon]$. The first step is to understand the composition of maps from $\mathrm{Spm}(L[\varepsilon])$ with the map to weight space.

We start by understanding maps to the 1 dimensional weight space, $\mathcal{W}$. By the functor of points definition of weight space, these maps correspond to characters $\tilde{\chi} : \mathbb{Z}_p^{\times} \to L[\varepsilon]^{\times}$, extending a fixed character $\chi : \mathbb{Z}_p^{\times} \to L^{\times}$.
\begin{lemma}
All $\tilde{\chi}$ are given by $\tilde{\chi}(z) = \chi(z) (1 + t \varepsilon \log z)$ for some $t \in L$
\end{lemma}
\begin{proof}
Dividing $\tilde{\chi}$ by $\chi$ gives another character, $1 + \varepsilon \psi$. By the group homomoprhism property, $\psi : \mathbb{Z}_p^{\times} \to L$ satisfies $\psi(ab) = \psi(a) + \psi(b)$. By considering the action of $\psi$ on a topological generator of $(1 + p\mathbb{Z}_p)^{\times}$, it follows $\psi|_{(1 + p\mathbb{Z}_p)^{\times}}$ is $t \log$ for some $t \in L$. Combining this with the action on roots of unity shows the two functions agree for all of $\mathbb{Z}_p^{\times}$. 
\end{proof}

Extending this to the product of two weight spaces shows that the $L[\varepsilon]$-points on the eigenvariety have weights $((1 + t_1\varepsilon \log)\chi_1, (1 + t_2\varepsilon \log)\chi_2)$, where $(\chi_1, \chi_2)$ is the weight of the original point, and $(t_1, t_2) \in \mathbb{Q}_p^2$. This weight will be denoted as $(\chi_1 + t_1\varepsilon, \chi_2 + t_2\varepsilon)$.

The $L[\varepsilon]$-points on the eigenvariety correspond to eigenvalue systems for $\mathbb{T}$ acting on $H^1_c(\mathbb{D}_{\chi_1 + t_1\varepsilon, \chi_2 + t_2\varepsilon}(L[\varepsilon]))$. We observe that the eigenvalues of a class being valued in $L[\varepsilon]$ is the same as the class being a generalised eigenvector of rank 2 over $L$, along with some compatibility with multiplication by $\varepsilon$.

As $L[\varepsilon]$ is a finite dimensional $L$-algebra, $\mathbb{D}_{\chi_1 + t_1\varepsilon, \chi_2 + t_2\varepsilon}(L[\varepsilon])$ is isomorphic to $\mathrm{Hom}_{L[\varepsilon]}(\mathcal{A}(L[\varepsilon]), L[\varepsilon]) \cong \mathrm{Hom}_L(\mathcal{A}(L), L[\varepsilon])$. As a result, there is a short exact sequence of $L$-vector spaces.
$$0 \to \mathbb{D}_{\chi_1, \chi_2}(L) \xrightarrow{\times \varepsilon} \mathbb{D}_{\chi_1 + t_1\varepsilon, \chi_2 + t_2\varepsilon}(L[\varepsilon]) \xrightarrow{\varepsilon = 0} \mathbb{D}_{\chi_1, \chi_2}(L) \to 0$$

\begin{lemma}
This short exact sequence is a short exact sequence of $\Sigma_0(p)^2$-modules.
\end{lemma}
\begin{proof}
Exactness follows from the statement for vector spaces, so it remains to check the module structure. It is enough to check for elements of the form $(\sigma, 1)$ and $(1, \sigma)$. These two cases are symmetric and so we only check $(\sigma, 1)$. For $\times \varepsilon$, take $\mu \in \mathbb{D}_{\chi_1, \chi_2}(L)$ and let $\sigma = \begin{psmallmatrix}a & b \\ c & d\end{psmallmatrix}
$.
$$(\varepsilon \times (\mu \cdot (\sigma, 1)))(f) = \varepsilon \mu\left(\chi_1(a + cx) f\left(\frac{dx + b}{a + cx}, y\right)\right)$$
$$= \mu\left(\varepsilon \chi_1(a + cx) f\left(\frac{dx + b}{a + cx}, y\right)\right)$$
Since $\varepsilon^2 = 0$, $\varepsilon \chi_1 = \varepsilon (1 + t_1\varepsilon \log)\chi_1$ and the following equalities hold.
$$(\varepsilon \times (\mu \cdot (\sigma, 1)) = \mu\left(\varepsilon \chi_1(a + cx) (1 + t_1\varepsilon \log(a + cx)) f\left(\frac{dx + b}{a + cx}, y\right)\right)$$
$$= (\varepsilon \mu)\left(\chi_1(a + cx) (1 + t_1\varepsilon \log(a + cx)) f\left(\frac{dx + b}{a + cx}, y\right)\right)$$
$$= ((\varepsilon \times \mu) \cdot (\sigma, 1))(f)$$
The statement for the map $\varepsilon = 0$ follows by definition, as $\chi_i + t_i\varepsilon$ give $\chi_i$ under this map.
\end{proof}

This short exact sequence gives a Hecke-equivariant long exact sequence of cohomology. The following subsection of the sequence will be useful.
$$H^1_c(\mathbb{D}_{\chi_1, \chi_2}(L)) \to H^1_c(\mathbb{D}_{\chi_1 + t_1\varepsilon, \chi_2 + t_2\varepsilon}(L[\varepsilon])) \to H^1_c(\mathbb{D}_{\chi_1, \chi_2}(L)) \xrightarrow{\delta} H^2_c(\mathbb{D}_{\chi_1, \chi_2}(L))$$
The first map is an injection, as $H^0_c(\mathbb{D}_{\chi_1, \chi_2}(L)) = 0$.

The eigenvariety parameterises eigenvalue systems appearing in these cohomology groups. The following proposition shows that deforming eigenvalues is the same thing as deforming classes of symbols. As with all other propositions in this section, the assumptions of this proposition are satisfied by Bianchi cusp forms of non-critical slope, since these have finite slope and appear in 1-dimensional eigenspaces in both degrees of cohomology.

\begin{proposition}
Let $\pi$ be a system of eigenvalues, with finite slope, such that $H^1_c(\mathbb{D}_{\chi_1, \chi_2}(L))_{\pi}$ is 1-dimensional, and $\phi$ a non-zero element. The system $\pi$ deforms along $(t_1, t_2) \neq (0, 0)$ if and only if $\phi$ is in the image of the map $H^1_c(\mathbb{D}_{\chi_1 + t_1\varepsilon, \chi_2 + t_2\varepsilon}(L[\varepsilon])) \to H^1_c(\mathbb{D}_{\chi_1, \chi_2}(L))$.
\end{proposition}
\begin{proof}
To find a deformation of the eigenvalue system, $H^1_c(\mathbb{D}_{\chi_1 + t_1\varepsilon, \chi_2 + t_2\varepsilon}(L[\varepsilon]))_{\pi}$ must be 2-dimensional over $L$. As the long exact sequence is Hecke-equivariant, it restricts to eigenspaces. In particular, the following is exact.
$$0 \to H^1_c(\mathbb{D}_{\chi_1, \chi_2}(L))_{\pi} \to H^1_c(\mathbb{D}_{\chi_1 + t_1\varepsilon, \chi_2 + t_2\varepsilon}(L[\varepsilon]))_{\pi} \to H^1_c(\mathbb{D}_{\chi_1, \chi_2}(L)))_{\pi}$$
By assumption, both the first and third terms are 1 dimensional. For the middle term to be 2 dimensional, it must surject onto the third term, and so $\phi$ is in the image. 

Conversely, suppose $\phi$ is in the image and let $\tilde{\phi}$ be a lift. Let $h$ be larger than the slope of $\tilde{\phi}$, then $H^1_c(\mathbb{D}_{\chi_1 + t_1\varepsilon, \chi_2 + t_2\varepsilon}(L[\varepsilon]))^{\leq h}$ is finite dimensional as an $L$-vector space. Combining this with the Hecke equivariance of the map shows we can assume $\tilde{\phi}$ is an eigensymbol. Take $a \in \mathbb{Z} \setminus \{\pm 1\}$ coprime to the level, and $b \in \mathcal{O}_L$ be such that $b$ and $\bar{b}$ are coprime to the level, and moreover, $c = \frac{\bar{b}}{b}$ is not a root of unity. The eigenvalue of $\langle a \rangle$ on $\tilde{\phi}$ is $\chi_1(a) \chi_2(a) (1 + (t_1 + t_2)\varepsilon\log a)$, and for $\langle c \rangle$, it is $\chi_1(c) \chi_2(\bar{c}) (1 + t_1\varepsilon \log c + t_2\varepsilon \log \bar{c}) = \chi_1(c) \chi_2(c)^{-1} (1 + (t_1 - t_2)\varepsilon \log c)$. This gives a deformation of the eigenvalue system, unless both $\chi_1(a) \chi_2(a) (t_1 + t_2)\varepsilon\log a = 0$ and $\chi_1(c) \chi_2(c)^{-1} (t_1 - t_2)\varepsilon\log c = 0$. As $a$ and $c$ are coprime to the level, $\chi_i(a)$, $\chi_i(c)$ are non-zero, as are $\log a$, $\log c$. As $(t_1, t_2) \neq (0, 0)$, at least one of $t_1 + t_2$ and $t_1 - t_2$ is non-zero, and the eigenvalue system $\pi$ deforms.
\end{proof}
If $t_1 = t_2 = 0$, then the short exact sequence of modules splits. In this case, the symbol deforms if and only if the eigenvalue system corresponds to a generalised eigensystem in $H^1_c(\mathbb{D}_{\chi_1, \chi_2}(L))$. If $\phi$ is a Hecke eigensymbol with $T_{\mathfrak{q}}$ eigenvalue, $a_{\mathfrak{q}}$, and $\psi$ is a generalised eigensymbol such that $(T_{\mathfrak{q}} - a_{\mathfrak{q}}) \psi$ is a multiple of $\phi$, then the correspondence sends $(\phi, \psi)$ to $\phi + \varepsilon \psi$. As Bianchi cusp forms which are not of critical slope have 1 dimensional eigenspaces, we can assume that $(t_1, t_2) \neq (0, 0)$ going forward. 

We use these methods to show a previously known result: the non-critical cuspidal components of the Bianchi eigenvariety contains no isolated points at cohomological weights \cite{hannew}.
\begin{theorem}
Let $\pi$ be a system of eigenvalues such that $H^1_c(\mathbb{D}_{\chi_1, \chi_2}(L))_{\pi}$ and $H^2_c(\mathbb{D}_{\chi_1, \chi_2}(L))_{\pi}$ are both 1 dimensional, then $\pi$ deforms along some non-zero direction in weight space, and so the tangent space to the eigenvariety has positive dimension.
\label{nopts}
\end{theorem}
For this, we need a lemma controlling the behaviour of the connecting homomorphism, $\delta : H^1_c(\mathbb{D}_{\chi_1, \chi_2}) \to H^2_c(\mathbb{D}_{\chi_1, \chi_2})$.
\begin{lemma}
Let $\Delta : H^1_c(\mathbb{D}_{\chi_1, \chi_2}) \times \mathbb{Q}_p^2 \to H^2_c(\mathbb{D}_{\chi_2, \chi_2})$ be such that $\Delta(\_, (t_1, t_2))$ is the boundary homomorphism associated to $(\chi_1 + t_1 \varepsilon, \chi_2 + t_2 \varepsilon)$, then $\Delta$ is linear in the second argument.
\label{lindef}
\end{lemma}
\begin{proof}
Let $L[\varepsilon_1, \varepsilon_2]$ denote $L[\varepsilon_1, \varepsilon_2] / (\varepsilon_1^2, \varepsilon_1 \varepsilon_2, \varepsilon_2^2)$ by analogy with $L[\varepsilon]$. We can similarly define a $\Sigma_0(p)^2$-module $\mathbb{D}_{\chi_1 + \varepsilon_1, \chi_2 + \varepsilon_2}(L[\varepsilon, \varepsilon_2])$, where the action is via the pair of characters $(\chi_1 + \varepsilon_1 \log, \chi_2 + \varepsilon_2 \log)$. As before, there is a short exact sequence of modules, where the first map sends $(\mu_1, \mu_2)$ to $\varepsilon_1 \mu_1 + \varepsilon_2 \mu_2$, and the second map sends $\varepsilon_1$ and $\varepsilon_2$ to 0.
$$0 \to \mathbb{D}_{\chi_1, \chi_2}(L)^2 \to \mathbb{D}_{\chi_1 + \varepsilon_1, \chi_2 + \varepsilon_2}(L[\varepsilon_1, \varepsilon_2]) \to \mathbb{D}_{\chi_1, \chi_2}(L) \to 0$$
This short exact sequence is compatible with the one corresponding to $(\chi_1 + t_1 \varepsilon, \chi_1 + t_2\varepsilon)$, in that the following commutative diagram is commutative, where the first vertical map is $(\mu_1, \mu_2) \mapsto (t_1 \mu_1 + t_2 \mu_2)$, the second is given by $\varepsilon_1 \mapsto t_1 \varepsilon$ and $\varepsilon_2 \mapsto t_2 \varepsilon$, and the third is the identity.
$$\begin{tikzcd}
   0 \arrow[r] & \mathbb{D}_{\chi_1, \chi_2}(L)^2 \arrow[r]\arrow[d] & \mathbb{D}_{\chi_1 + \varepsilon_1, \chi_2 + \varepsilon_2}(L[\varepsilon_1, \varepsilon_2]) \arrow[r]\arrow[d] & \mathbb{D}_{\chi_1, \chi_2}(L) \arrow[r]\arrow[d] & 0 \\
   0 \arrow[r] & \mathbb{D}_{\chi_1, \chi_2}(L) \arrow[r] & \mathbb{D}_{\chi_1 + t_1\varepsilon, \chi_2 + t_2\varepsilon}(L[\varepsilon]) \arrow[r] & \mathbb{D}_{\chi_1, \chi_2}(L) \arrow[r] & 0 \\
  \end{tikzcd}
$$
This induces the following commutative square of cohomology.
$$\begin{tikzcd}
   H^1_c(\mathbb{D}_{\chi_1, \chi_2}(L)) \arrow[d]\arrow[r, "\delta'"] & H^2_c(\mathbb{D}_{\chi_1, \chi_2}(L)^2) \arrow[d] \\
   H^1_c(\mathbb{D}_{\chi_1, \chi_2}(L)) \arrow[r, "\delta"] & H^2_c(\mathbb{D}_{\chi_1, \chi_2}(L))\\
  \end{tikzcd}
$$
The top right corner of this square is isomorphic to $H^2_c(\mathbb{D}_{\chi_1, \chi_2}(L))^2$, and so with this identification, and using that the lefthand vertical map is an equality, $\Delta$ can be recognised as the composition of $\delta'$ and $(\phi_1, \phi_2) \mapsto t_1 \phi_1 + t_2 \phi_2$. The only dependence on $(t_1, t_2)$ comes from this second map, which is linear in $(t_1, t_2)$, and so the result follows.
\end{proof}

We can now use this to prove the theorem.
\begin{proof}[Proof of Theorem~\ref{nopts}]
Let $\phi$ be an eigensymbol in the eigenvalue system $\pi$. By the proposition and following remark the eigenvalue system will only deform if $\phi$ is in the image of $H^1_c(\mathbb{D}_{\chi_1 + t_1\varepsilon, \chi_2 + t_2\varepsilon}(L[\varepsilon])) \to H^1_c(\mathbb{D}_{\chi_1, \chi_2}(L))$ for some choice of $(t_1, t_2) \neq (0, 0)$. From the long exact sequence, this happens precisely if $\Delta(\phi, (t_1, t_2)) = 0$. As $\Delta(\phi, (t_1, t_2))$ is contained in $H^2_c(\mathbb{D}_{\chi_1, \chi_2}(L))_{\pi}$ by Hecke equivariance, $\Delta(\phi, \_)$ gives a map from the tangent space to $L$. The preceeding lemma shows this is a linear map, and as the tangent space of the weight space is 2-dimensional, it has a non-trivial kernel. 
\end{proof}

For computational purposes, having to compute all moments of the distribution in the deformation part is unfeasible. We show that it is enough to compute only a finite number of the moments in a manner analogous to the classical projection of an overconvergent modular symbol. We now assume that $\chi_1(x) = \chi_2(x) = x^k$ for some non-negative integer $k$ as Bianchi cusp forms are restricted to parallel weights. Define $\mathbb{D}^0_{k, k}(L) = \{\mu \in \mathbb{D}(L) \mid \mu(x^i y^j) = 0, 0 \leq i, j \leq k\}$. This subspace is $\Sigma_0(p)^2$-invariant, therefore $\varepsilon \mathbb{D}^0_{k, \ell}(L) \subset \mathbb{D}_{k + t_1\varepsilon, k + t_2\varepsilon}(L[\varepsilon])$ is also $\Sigma_0(p)^2$-invariant by the $\Sigma_0(p)^2$-invariance of $\times \varepsilon$. The quotient by this subspace will be denoted $\mathbb{D}_{k+t_1\varepsilon, k+t_2\varepsilon}(L[\varepsilon])'$. This space fits into a similar short exact sequence as before, which is again, $\Sigma_0(p)^2$-invariant. 
$$0 \to V_{k, k}(L) \xrightarrow{\times \varepsilon} \mathbb{D}_{k + t_1\varepsilon, k + t_2\varepsilon}(L[\varepsilon])' \xrightarrow{\varepsilon = 0} \mathbb{D}_{k, k}(L) \to 0$$

We now show that this is enough to recover all the information we are interested in.
\begin{theorem}
Let $\phi$ be a $U_p$-eigensymbol in $H^1_c(\mathbb{D}_{k, k}(L))$ of slope strictly less than $k + 1$, then $\phi$ deforms in the direction $(t_1, t_2)$ if and only if it is in the image of the map $H^1_c(\mathbb{D}_{k + t_1\varepsilon, k + t_2\varepsilon}(L[\varepsilon])') \to H^1_c(\mathbb{D}_{k, k}(L))$.
\end{theorem}
\begin{proof}
Combining the two short exact sequences gives a commutative diagram, where the vertical maps are the quotients discussed previously.
$$\begin{tikzcd}
   0 \arrow[r] & \mathbb{D}_{k, k}(L) \arrow[r]\arrow[d] & \mathbb{D}_{k + t_1\varepsilon, k + t_2\varepsilon}(L[\varepsilon]) \arrow[r]\arrow[d] & \mathbb{D}_{k, k}(L) \arrow[r]\arrow[d] & 0 \\
   0 \arrow[r] & V_{k, k}(L) \arrow[r] & \mathbb{D}_{k + t_1\varepsilon, k + t_2\varepsilon}(L[\varepsilon])' \arrow[r] & \mathbb{D}_{k, k}(L) \arrow[r] & 0 \\
  \end{tikzcd}
$$
As this is commutative, by the naturality of long exact sequences, the following diagram is also commutative. The vertical maps are maps induced by the quotients.
$$\begin{tikzcd}
   H^1_c(\mathbb{D}_{k + t_1\varepsilon, k + t_2\varepsilon}(L[\varepsilon])) \arrow[r]\arrow[d] & H^1_c(\mathbb{D}_{k, k}(L)) \arrow[d]\arrow[r, "\delta"] & H^2_c(\mathbb{D}_{k, k}) \arrow[d] \\
   H^1_c(\mathbb{D}_{k + t_1\varepsilon, k + t_2\varepsilon}(L[\varepsilon])') \arrow[r] & H^1_c(\mathbb{D}_{k, k}(L)) \arrow[r, "\delta'"] & H^2_c(V_{k, k}(L))\\
  \end{tikzcd}
$$

The symbol $\phi$ deforms when it is in the image of the first map in the top row. By commutativity of the diagram, if $\phi$ deforms, then it is in the image of the first map in the second row.

Assume now that it is in the image of the first map in the second row. As $\phi$ has slope less than $k + 1$, so does $\delta(\phi)$, and the right hand map is an isomorphism when restricted to the subspace spanned by $\delta(\phi)$ (by the control theorem in \cite{williams}). By exactness, $\delta'(\phi) = 0$, therefore $\delta(\phi) = 0$. This shows $\phi$ deforms.
\end{proof}
As the quotient map on distribution spaces is $\Sigma_0(p)^2$-invariant, the induced map $H^1_c(\mathbb{D}_{k + t_1\varepsilon, k + t_2\varepsilon}(L[\varepsilon])) \to H^1_c(\mathbb{D}_{k + t_1\varepsilon, k + t_2\varepsilon}(L[\varepsilon])')$ is Hecke equivariant. This means that deformations to eigenvalues can be computed in the quotient space as well as the direction. 

\section{Computational Results}
The previous section suggests computing deformation directions for eigenvalue systems by the following method. We use modular symbols to make the cohomology spaces explicit. 
\begin{enumerate}
 \item Compute a generating set and the Manin relations for $\mathrm{Div}^0(\mathbb{P}^1(K))$.
 \item Solve the Manin relations in the case of $V_{k, k}$ to compute $H^1_c(V_{k, k})$ as $\mathrm{Symb}_{\Omega}(V_{k, k})$.
 \item Identify a modular symbol with the specified eigenvalue system.
 \item Compute an overconvergent lift of the modular symbol, that is, a lift to $\mathrm{Symb}_{\Omega}(\mathbb{D}_{k, k}(L))$.
 \item As in the definition of the connecting homomorphism, determine the failure of the trivial deformation to satisfy the Manin relations when deforming along $(1, 0)$ and $(0, 1)$.
 \item Find which combination of the errors computed in the previous step can be accounted for using the Manin relations for $V_{k, k}$.
\end{enumerate}

The first step has been done in for Euclidean fields \cite{cremona}. Subsequent work has extended this to the remaining class number 1 fields \cite{whitley}, and also to higher class numbers \cite{bygott}. Let $\gamma_1, \ldots, \gamma_n$ be a generating set for $\mathrm{Div}^0(\mathbb{P}^1(K))$, then there is a linear map $R : \left(K\gamma_1 + \ldots + K \gamma_n\right)\otimes V_{k, k} \to K^m$ encoding the Manin relations, where $m$ is the number of relations. Solving the relations consists of finding the kernel of this map. Using a database of Bianchi modular forms, such as the LMFDB \cite{lmfdb}, it is straightforward to find a Hecke operator where the eigenvalue is unique. The modular symbol is then the kernel of an operator. As steps 1-3 are all linear algebra (over a number field), they can be carried out exactly on a computer. 

Overconvergent lifts can be computed using finite approximation modules and repeated application of the $U_p = U_{\mathfrak{p}} U_{\bar{\mathfrak{p}}}$ operator and ``dividing'' by the eigenvalue \cite{williams}. Unfortunately, this is an infinite process. 

\begin{proposition}
Steps 4-5 can be computed to finite precision, and the error arising is bounded.
\end{proposition}
\begin{proof}
Truncating the process at the $N$th finite approximation module computes the values of $\mu(x^i y^j) \pmod{\pi_L^{N - i - j}}$, for a uniformiser, $\pi_L$, of $\mathcal{O}_L$, since the lifting process respects the maps between finite approximation modules \cite{williams}. 

Computing the failure of the trivial deformation to satisfy the Manin relations requires evaluating the distributions attached to the generators on power series of the form $\log(a + pcx) x^i y^j$ and $\log(a + pcy)x^i y^j$ for $0 \leq i, j \leq k$. The coefficients arising in these power series are of the form $p^n c^n / n$, which have valuation at least $n - \nu_p(n)$. For $n + i + j \geq N$, these contributions have valuation at least $n - \nu_p(n)$, and if $n + i + j < N$, the contributions have a rounding error with valuation at least $(N - i - j)  - \nu_p(n)$. These errors can therefore be made smaller than a given bound by choosing a suitable $N$. 
\end{proof}

Finally, to determine which direction the symbol deforms in, we use that the failure to deform is linear in direction, c.f. Lemma~\ref{lindef}. A solution to the Manin relations will exist if there is a set of terms (valued in $\varepsilon V_{k, k}$) that can be added to the trivial deformation so that the Manin relations are satisfied. In particular, to deform along the direction $(1, t)$, the equation $R x = E_1 + t E_2$ must have a solution, where the $E_i$ are the failures to deform along $(1, 0)$ and $(0, 1)$ respectively. This is now linear algebra and can be solved numerically. As the right hand side is only known modulo a power of $p$, the most reliable approach is to use algorithms for solving linear systems over $\mathbb{Z}/p^n \mathbb{Z}$.

We implemented the above for the Euclidean quadratic imaginary fields, using SageMath \cite{sagemath}. The source code is available at \url{https://github.com/jameswrawson/OCBianchiSymbols}. Many steps of the process are simplified in this case, such as computing Hecke operators, as well as it only being necessary to work with a single space of modular symbols, rather than an ensemble.

For a modular form of level $\Gamma_0(\mathcal{N})$, weight $k$ and working precision $d$, we can estimate the complexity of this process. 
\begin{theorem}
Steps 1-6 can be computed with a time complexity polynomial in $N = \mathrm{Nm}(\mathcal{N})$, $k$, $d$ and $p$.
\end{theorem}
\begin{proof}
The index of $\Gamma_0(\mathcal{N})$ inside $\mathrm{SL}_2(\mathcal{O}_K)$ is slightly larger than $N$, it is $N \prod_{\mathfrak{q} | \mathcal{N}} \left(1 + \frac{1}{\mathrm{Nm}(\mathfrak{q})}\right)$, where the product is over prime ideals. The product is at most $\prod_{q | N} (1 + \frac{1}{q})^2$, which is bounded by $N^{\delta}$ for any $\delta > 0$. The number of generators and relations needed is proportional to this index, as each coset is attached to a generator and a collection of relations. However, to compute the Manin relations for $V_{k, k}$, $(k + 1)^2$ copies of each relation are needed. The resulting matrix is therefore approximately square of size $k^2 N^{1 + \delta}$. The complexity of computing the kernel of such a matrix by Gaussian elimination is then $\mathcal{O}(k^6 N^{3 + 3 \delta})$, with the possibility of more efficient algorithms being available to utilise that this matrix is sparse. 

The complexity of selecting a modular symbol from the space by Hecke operators is hard to quantify, as it depends on the dimension of the cohomology space and on the Hecke operator used.
For each matrix comprising the Hecke operator $T_{\mathfrak{q}}$, the value of the modular symbol must be determined on the image of every generator under this matrix. The resulting paths must then be split into unimodular paths via continued fractions, typically giving $\log(qN)$ paths, where $q$ is the norm of $\mathfrak{q}$. On these paths, the value can be computed by $\Gamma_0(\mathcal{N})$-invariance. This represents $\mathcal{O}(q\log(qN)N^{1 + \delta})$ operations in all. Since the dimension of the space of modular symbols is at most $\mathcal{O}(k^2 N^{1 + \delta})$, computing the Hecke operator on the whole is space can still be done in polynomial time.

Computing the $U_{\mathfrak{p}}$ operator to compute the lifts is similar, except that computing the values on the unimodular paths is less straightforward. Each path needs $\mathcal{O}(d^2)$ values calculating (the stored moments), each needing the first $\mathcal{O}(d^2)$ terms of the power series of a rational function. This represents a total of around $pd^4N^{1 + \delta} \log(pN)$ operations. This is repeated again for $U_{\bar{\mathfrak{p}}}$, and the pair is then repeated $d$ times. In total, this represents $\mathcal{O}(p d^5 N^{1 + 2\delta})$ operations.

Determining the errors arising along the cardinal directions requires computing the error for each relation with values in $V_{k, k}$. For each relation, the first $kd$ terms of a power series are needed (either $(a + cx)^k (\bar{a} + \bar{c}y)^k \log(1 + pax)$ or its mirror image). This gives a total of $\mathcal{O}(k^3 d N^{1 + \delta})$ operations. Finally, solving the linear equation for the direction is comparable to the first stages of the calculation. 
\end{proof}

In practice, the most time consuming step is 4, computing the overconvergent lift. This appears to be a result of the sensitivity to the precision (growing like $d^5$) as well as a comparatively large implicit constant. Although the complexity of step 3 is hard to give explicitly, in the examples given below, it was no more intensive than any other step.

The results of this process in a few cases in weight (2, 2) (equivalently, $k = 0$), are illustrated in Figure~\ref{bianchidef}. In these cases the primes used do not divide the level, so the ordinary $p$-stabilisation is used. As expected, all of these have a definite deformation direction, and none deform along the parallel weight direction, verifying the eigenvariety is smooth at these points and that the Calegari-Mazur conjecture holds for these families.

\begin{figure}
\begin{tabular}{|p{2cm}|p{2cm}|p{2cm}|p{4cm}|}
\hline
Field & Level & Prime & Direction of deformation $(1, \_)$\\\hline
$\mathbb{Q}(\sqrt{-11})$ & $(9 + \sqrt{-11})$ \newline norm 92 & $\left(\frac{1 + \sqrt{-11}}{2}\right)$\newline norm 3 & $2 + 2\times3 + 3^3 + 2\times3^4 + \mathcal{O}(3^5)$ \\\hline
$\mathbb{Q}(\sqrt{-11})$ & $(6 - \sqrt{-11})$\newline norm 47 & $\left(\frac{-3 + \sqrt{-11}}{2}\right)$\newline norm 5 & $3 + 2\times5 + 4\times5^2 + 5^3 + \mathcal{O}(5^4)$ \\\hline
$\mathbb{Q}(\sqrt{-11})$ & $\left(\frac{9 + 5\sqrt{-11}}{2}\right)$\newline norm 89 & $\left(\frac{3 + \sqrt{-11}}{2}\right)$\newline norm 5 & $4 + 4\times5 + \mathcal{O}(5^4)$ \\\hline
$\mathbb{Q}(\sqrt{-3})$ & $\left(\frac{17 + \sqrt{-3}}{2}\right)$\newline norm 73 & $\left(\frac{5 + \sqrt{-3}}{2}\right)$\newline norm 7 & $2 + 2\times7 + 3 \times 7^2 + 4\times7^3 + \mathcal{O}(7^4)$ \\\hline
\end{tabular}
\caption{Deformation directions for some Bianchi cusp forms of weight (2, 2).}
\label{bianchidef}
\end{figure}
\section{An Analogue for Elliptic Modular Forms}
The question of tangent spaces to the Coleman-Mazur eigencurve is also of interest, not for the direction over the (one dimensional) weight space, but for the derivatives of the Hecke eigenvalues. In particular, the $p$-adic derivative, with respect to the weight, of the $U_p$-eigenvalue is closely related to the $\mathcal{L}$-invariant of the modular form, $\mathcal{L} = -2 p^{-k / 2} a_p'$ for a form of weight $k + 2$. For cusp forms corresponding to elliptic curves with split multiplicative reduction at $p$, there is an interpretation of $\mathcal{L}$ in terms of the Tate parameterisation of the elliptic curve. More importantly, $\mathcal{L}$-invariants relate derivatives of $p$-adic $L$-functions and classical $L$-values \cite{stevens}. In Urban's proof of the Iwasawa Main Conjecture for adjoint representations \cite{urban}, the non-vanishing of $\mathcal{L}$-invariants for modular forms corresponding to elliptic curves is used to conclude equality in that case. We will use our tangent space methods to compute $\mathcal{L}$-invariants.

By direct analogy with the Bianchi case, computing deformations of eigenvalue systems consists of computing deformations to the modular symbols. In this case, there are two simplifications: there is no need to determine a deformation direction; and simpler generating sets and relations exist, which can be computed from fundamental domains \cite{pollstev}. 

We separately implemented this algorithm for elliptic modular forms. The results in a few cases are summarised below. Notably, in each example, the $\mathcal{L}$-invariant was non-zero. We start by tabulating these in Figure~\ref{linvars}, before also giving the deformations of other coefficients, whose properties are less well understood. The first two values of Figure~\ref{linvars} can be compared against the values computed directly from the elliptic curves (after scaling by the factor of $-\frac{1}{2}$). The third value represents an elliptic curve with non-split multiplicative reduction. Using the formula for elliptic curves with split multiplicative reduction \cite{stevens}, $\log(q)/\nu_p(q)$, where $q$ is the Tate parameter, gives the negative of the $\mathcal{L}$-invariant for this example, which can be seen by twisting by a suitable quadratic character. The last example is the ordinary 3-stabilisation of the cusp form specified, this is due to a limitation of the software used to compute overconvergent lifts of modular symbols. Figures~\ref{def11},~\ref{def55},~\ref{def15}, and~\ref{def5} show the remaining Hecke eigenvalues computed.

\begin{figure}
\begin{tabular}{|p{1.4cm}|p{1.4cm}|p{1.4cm}|p{1.4cm}|p{4.2cm}|}
\hline
Weight & Level & LMFDB identifier & Prime & Derivative of the $U_p$ eigenvalue \\\hline
2 & 11 & 11.2.a.a & 11 & $8\times11 + 2 \times 11^2 + 7 \times 11^3 + 11^4 + 7 \times 11^5 + 11^6 + 3 \times 11^7 + \mathcal{O}(11^8)$ \\\hline
2 & 55 & 55.2.a.a & 5 & $4 \times 5 + 2 \times 5^3 + 4 \times 5^4 + 2 \times 5^5 + 4 \times 5^7 + \mathcal{O}(5^8)$ \\\hline
2 & 15 & 15.2.a.a & 3 & $2\times3 + 3^2 + 2 \times 3^3 + 2 \times 3^5 + 3^6 + \mathcal{O}(3^7)$ \\\hline
4 & 5 & 5.4.a.a & 3 & $2 \times 3 + 2 \times 3^2 + 2 \times 3^5 + 3^6 + \mathcal{O}(3^7)$ \\\hline
\end{tabular}
\caption{Derivatives of the $U_p$ operator for a selection of cusp forms}
\label{linvars}
\end{figure}

\begin{figure}
\begin{tabular}{|p{2cm}|p{9cm}|}
\hline
Operator & Deformed Eigenvalue\\\hline
$T_2$ & $2\times11 + 11^2 + 2\times 11^3 + 8 \times 11^4 + 10 \times 11^5 + 7 \times 11^6 + 9 \times 11^7 + \mathcal{O}(11^8)$ \\\hline
$T_3$ & $10 \times 11 + 6 \times 11^2 + 11^3 + 10 \times 11^4 + 7 \times 11^5 + 8 \times 11^6 + 2 \times 11^7 + \mathcal{O}(11^8)$ \\\hline
$T_5$ & $11 + 10 \times 11^2 + 11^4 + 9 \times 11^5 + 6 \times 11^6 + 3 \times 11^7 + \mathcal{O}(11^8)$ \\\hline
$T_7$ & $4\times11 + 4 \times 11^3 + 11^4 + 3 \times 11^5 + 5 \times 11^6 + 5 \times 11^7 + \mathcal{O}(11^8)$ \\\hline
$T_{13}$ & $8 \times 11 + 9 \times 11^2 + 11^4 + 5 \times 11^5 + 6 \times 11^6 + 8 \times 11^7 + \mathcal{O}(11^8)$ \\\hline
\end{tabular}
\caption{Derivatives of some Hecke operators for the cusp form with LMFDB label 11.2.a.a}
\label{def11}
\end{figure}

\begin{figure}
\begin{tabular}{|c|c|}
\hline
Operator & Deformed Eigenvalue\\\hline
$T_2$ & $5 + 3 \times 5^2 + 4 \times 5^3 + 5^4 + 2 \times 5^5 + 5^7 + \mathcal{O}(5^8)$ \\\hline
$T_3$ & $3 \times 5 + 4 \times 5^2 + 2 \times 5^3 + 4 \times 5^5 + 3 \times 5^7 + \mathcal{O}(5^8)$ \\\hline
$T_7$ & $3 \times 5 + 4 \times 5^2 + 2 \times 5^3 + 3 \times 5^5 + 2 \times 5^6 + 5^7 + \mathcal{O}(5^8)$ \\\hline
$U_{11}$ & $4 \times 5 + 2 \times 5^3 + 2 \times 5^4 + 3 \times 5^6 + \mathcal{O}(5^8)$ \\\hline
$T_{13}$ & $2 \times 5^2 + 4 \times 5^3 + 3 \times 5^4 + 5^7 + \mathcal{O}(5^8)$ \\\hline
\end{tabular}
\caption{Derivatives of some Hecke operators for the cusp form with LMFDB label 55.2.a.a}
\label{def55}
\end{figure}

\begin{figure}
\begin{tabular}{|c|c|}
\hline
Operator & Deformed Eigenvalue\\\hline
$T_2$ & $3 + 3^2 + 2 \times 3^3 + 3^4 + 3^5 + \mathcal{O}(3^7)$ \\\hline
$U_5$ & $2 \times 3 + 3^3 + 2 \times 3^5 + 3^6 + \mathcal{O}(3^7)$ \\\hline
$T_7$ & $3 + 2 \times 3^2 + 3^4 + 2 \times 3^5 + 3^6 + \mathcal{O}(3^7)$ \\\hline
$T_{11}$ & $2 \times 3^2 + 2 \times 3^3 + 2 \times 3^4 + 2 \times 3^6 \mathcal{O}(3^7)$ \\\hline
$T_{13}$ & $2 \times 3^5 + \mathcal{O}(3^7)$ \\\hline
\end{tabular}
\caption{Derivatives of some Hecke operators for the cusp form with LMFDB label 15.2.a.a}
\label{def15}
\end{figure}

\begin{figure}
\begin{tabular}{|c|c|}
\hline
Operator & Deformed Eigenvalue\\\hline
$T_2$ & $3 + 2 \times 3^5 + \mathcal{O}(3^7)$ \\\hline
$U_3$ & $2 \times 3 + 2 \times 3^2 + 2 \times 3^5 + 3^6 + \mathcal{O}(3^7)$ \\\hline
$T_7$ & $3 + 3^3 + 3^4 +  2 \times 3^6 + \mathcal{O}(3^7)$ \\\hline
$T_{11}$ & $2 \times 3^2 + 2 \times 3^4 + 3^5 + 2 \times 3^6 + \mathcal{O}(3^7)$ \\\hline
$T_{13}$ & $2 \times 3^2 + 3^3 + 2 \times 3^4 + 3^5 + \mathcal{O}(3^7)$ \\\hline
\end{tabular}
\caption{Derivatives of some Hecke operators for the cusp form with LMFDB label 5.4.a.a}
\label{def5}
\end{figure}

\subsection*{Acknowledgements}
We would like to thank David Loeffler for suggesting this problem, as well as for his advice and oversight during this project. 

\printbibliography
\end{document}